\documentclass[11pt]{article}
\usepackage{t1enc}
\usepackage[latin1]{inputenc}
\usepackage[english]{babel}
\usepackage{latexsym}
\usepackage{amssymb}
\usepackage{euscript}
\begin{document}

\setlength{\arraycolsep}{.136889em}
\renewcommand{\theequation}{\thesection.\arabic{equation}}
\newtheorem{thm}{Theorem}[section]
\newtheorem{propo}{Proposition}[section]
\newtheorem{lemma}{Lemma}[section]
\newtheorem{corollary}{Corollary}[section]
\newtheorem{remark}{Remark}[section]
\def\begg{\begin{equation}}
\def\endd{\end{equation}}
\def\ep{\varepsilon}
\def\noo{n\to\infty}
\def\al{\alpha}
\def\be{\bf E}
\def\bp{\bf P}
\medskip
\centerline{\Large\bf Random walks on the two-dimensional K-comb lattice}

\bigskip\bigskip

\bigskip\bigskip

\renewcommand{\thefootnote}{1}

\noindent 
{\textbf{Endre Cs\'aki}}

\noindent
Alfr\'ed R\'enyi Institute of Mathematics,  Budapest, P.O.B. 127, H-1364, Hungary. E-mail address:
csaki.endre@renyi.hu

\bigskip
\noindent

\bigskip

\noindent {\textbf{Ant\'{o}nia F\"{o}ldes}
\newline
Department of Mathematics, College of Staten Island, CUNY, 2800
Victory Blvd., Staten Island, New York 10314, U.S.A.  E-mail
address: Antonia.Foldes@csi.cuny.edu

\bigskip

\centerline{\bf Abstract}

\medskip
\noindent We study the path behavior of the
 symmetric  walk  on some special comb-type subsets of ${\mathbb Z}^2$
which are obtained from ${\mathbb Z}^2$ by generalizing the comb having finitely many  horizontal lines instead of one.

\medskip


\medskip

 \vspace{.1cm}

\section{Introduction}

\renewcommand{\thesection}{\arabic{section}} \setcounter{equation}{0}

\setcounter{thm}{0} \setcounter{lemma}{0}

The anisotropic random walk has a huge literature. Some important early work on this topic is due to Heyde \cite{H},
\cite{H93}. In our papers \cite{CCFR13}, \cite {CFR15} we give an account of some of the relevant 
literature. An anisotropic walk is defined  as a nearest neighbor random walk 
on the square lattice $\mathbb{Z}^2$ of the plane with possibly unequal 
symmetric horizontal and vertical step probabilities, so that these 
probabilities depend only on the value of the vertical coordinate. 

A very important special case is the simple random walk on the 2-dimensional 
comb lattice which is obtained from ${\mathbb Z}^2$ by removing all horizontal 
lines off the $x$-axis.
More formally, consider the random 
walk $\{{\bf C}(N)=\left(C_1(N),C_2(N)\right);\, N=0,1,2,\ldots\}$
on $\mathbb{Z}^2$ with the transition probabilities for  $(k,j)\in{\mathbb Z^2}$, $N=0,1,2,\ldots.$ 

\begin{eqnarray} 
{\bf P}({\bf C}(N+1)&=&(k, j\pm1)|{\bf C}(N)=(k,j))=\frac{1}{2} \,\, {\rm if } \,\, j\neq 0  \nonumber \\
{\bf P}({\bf C}(N+1)&=&(k,\pm1)|{\bf C}(N)=(k,0))= 
{\bf P} ({\bf C}(N+1)=(k\pm1, 0)|{\bf C}(N)=(k,0))=\frac{1}{4} 
\end{eqnarray} 
 For a recent review of some related literature concerning this
simple random walk we refer to Bertacchi \cite{BE} and Cs\'aki {\it et al.} 
\cite{CCFR08}. In the latter paper we established a simultaneous strong
approximation for the two coordinates  of the random walk 
${\bf C}(N)=(C_1(N),C_2(N))$ that reads as follows.

\medskip\noindent
{\bf Theorem A} (\cite{CCFR08}) {\it On an appropriate probability space 
for the simple random walk 
\newline $\{{\bf C}(N)=(C_1(N),C_2(N)); N=0,1,2,\ldots\}$ on the 
two-dimensional comb lattice ${\mathbb C}^2,$ one can construct two 
independent standard Wiener processes $\{W_1(t);\, t\geq 0\}$, 
$\{W_2(t);\, t\geq 0\}$ so that, as $N\to\infty$, we have with any 
$\varepsilon>0$
$$
N^{-1/4}|C_1(N)-W_1(\eta_2(0,N))|+N^{-1/2}|C_2(N)-W_2(N)|
=O(N^{-1/8+\varepsilon})\quad a.s.,
$$
where $\eta_2(0,\cdot)$ is the local time process at zero of}
$W_2(\cdot)$.

In this paper we want to generalize this theorem,  by permitting to have finitely  many horizontal lines, instead of one,
 and also relax  the requirements  about  the probabilities along these lines as follows. 
 At first we select a finite set of permitted horizontal lines.
 We call this set 
$B$ containing the horizontal lines  at $y=m_j, j=1,2,\ldots, K.$
Then we consider the random  walk 
on $\mathbb{Z}^2$ with the transition probabilities  for $(k,j)\in{\mathbb Z^2}$, $N=0,1,2,\ldots.$

$${\bf P}({\bf C}(N+1)=(k, \ell\pm1)|{\bf C}(N)=(k,\ell))=\frac{1}{2} \,\, {\rm if } \,\, \ell \notin  B$$
and for   $m_j\in  B$ 
\begin{eqnarray}
{\bf P}({\bf C}(N+1)&=&(k,m_j\pm1)|{\bf C}(N)=(k,m_j))= p_j   \nonumber \\
 {\bf P} ({\bf C}(N+1)&=&(k\pm1,m_j)|{\bf C}(N)=(k,m_j))=\frac{1}{2}-p_j 
\end{eqnarray} 
 where $0< p_j < \frac{1}{2}. $

\medskip\noindent
\medskip 
 Unless otherwise stated, we assume also that 
${\mathbf C}(0)=(0,0)$. We will call the above random walk  as K-comb walk on ${\mathbb C}_K^2.$  Introduce 
the notations  $\alpha_j=2p_j, \,\, j= 1,2\ldots, K$ and  

\begg A_K= \sum_{j=1}^ K \frac{1-\alpha_j}{\alpha_j}.  \label{alfa}
\endd
\noindent
In our paper  \cite{CF20}   we considered more general   (much bigger) $B$  sets but only under the condition of $p_j=1/4$ for all $j$-s.
The main result of this paper is the following generalization of Theorem A:

 \noindent
{\begin{thm}  On an appropriate probability space 
for the  random walk \newline
 $\{{\bf C}(N)=(C_1(N),C_2(N)); N=0,1,2,\ldots\}$ on the 
two-dimensional K-comb lattice ${\mathbb C}_K^2,$ one can construct two 
independent standard Wiener processes $\{W_1(t);\, t\geq 0\}$, 
$\{W_2(t);\, t\geq 0\}$ so that, as $N\to\infty$, we have with any 
$\varepsilon>0$
$$
N^{-1/4}|C_1(N)-W_1(A_K\eta_2(0,N))|+N^{-1/2}|C_2(N)-W_2(N)|
=O(N^{-1/8+\varepsilon})\quad a.s.,
$$
where $\eta_2(0,\cdot)$ is the local time process at zero of $W_2(\cdot).$ \end{thm}

\medskip \noindent
{\bf Remark 1.1 } Observe that the above result does not depend on the positions of the $K$ horizontal lines, only on their numbers and $A_K$ defined in (\ref{alfa}).

\medskip \noindent
{\bf Remark 1.2 } In case $K=1$ and $p_1=1/4$ our Theorem coincides with Theorem A.
\medskip

The structure of this paper from now on is as follows. In Section 2
we give preliminary facts and results. In  Section 3, first we redefine the 
walk on ${\mathbb C}_K^2$ in terms of two independent simple symmetric walks, and prove our Theorem. 
In Section 4 some consequences will be discussed. 

\section{Preliminaries}
\renewcommand{\thesection}{\arabic{section}} \setcounter{equation}{0}

\setcounter{thm}{0} \setcounter{lemma}{0}
In this section we list some well-known results, and some new ones which 
will be used in the rest of the paper. In case of the known ones  we won't 
give the most general form of the results, just as much as we intend to use,
while the exact reference will also be provided for the interested
reader.

Let $\{X_i\}_{i\geq1}$ be a sequence of independent i.i.d. random variables, 
with ${\mathbf P}(X_i=\pm1)=1/2.$ Then the simple symmetric random walk on 
the line is defined as $S(n)=\sum_{i=1}^n X_i, $  and its local time is  
$\xi(j,n)=\#\{k: 0< k \leq n,  S(k)=j\}, \,\,n=1,2...,$  for any integer $j.$

For $\xi(n)=\sup_{x}\xi(x,n)$ we have Kesten's LIL for local time.

\medskip\noindent {\bf Lemma A} (Kesten \cite{K}) {\it For the maximal local 
time we have}
$$ \limsup_{n\to\infty} \frac{\xi(n)}{(2n\log \log n)^{1/2}}=
1\quad a.s.$$

Let $(W(t),\, t\geq 0)$ be a standard Wiener process (called also standard
Brownian motion). Its local time $(\eta(x,t),\, x\in \mathbb R,\, \, t\geq 0)$
(called Wiener local time or Brownian local time) is defined as
$$
\eta(x,t)=\lim_{\varepsilon\to 0}\frac{1}{2\varepsilon}
\int_0^t I\{W(s)\in (x-\varepsilon,x+\varepsilon)\}\, ds,
$$
where $I\{\cdot\}$ denotes the indicator function.

Concerning the increments of the Brownian motion, Brownian local time and their   random walk 
counterparts we quote the following
result from Cs\"org\H o and R\'ev\'esz \cite{CSR79}, (see in   \cite{RE} page 69),  Cs\'aki  and F\"oldes \cite{CF83}
and Cs\'aki {\it et al.}  \cite{CCFR83}.

\medskip\noindent 
{\bf Lemma B}  {\it Let $0<a_T\leq T$ be a non-decreasing 
function of $T$. Then, as $T\to\infty$, we have
$$
\sup_{0\leq t\leq T-a_T}\sup_{0\leq s\leq a_T}|W(t+s)-W(t)|=
O(a_T(\log(T/a_T)+\log\log T))^{1/2} \qquad a.s.
$$
$$
\sup_{0\leq t\leq T-a_T}|\eta(0,t+a_T)-\eta(0,t)|=
O(a_T(\log(T/a_T)+\log\log T))^{1/2} \qquad a.s.
$$
$$ \sup_{0\leq t\leq T-a_T}|\xi(0,t+a_T)-\xi(0,t)|=
O(a_T(\log(T/a_T)+\log\log T))^{1/2} \qquad a.s. ,
$$
where in the last line $t$ and  $T$ should be  integers.}

We quote the following simultaneous strong approximation result from 
R\'ev\'esz \cite{RE81}. 

\medskip\noindent
{\bf Lemma C} (\cite{RE81}) {\it On an appropriate probability space for a 
simple symmetric random walk 

\noindent
$\{S(n);\, n=0,1,2,\ldots\}$ with local time
$\{\xi(x,n);\, x=0,\pm1,\pm2,\ldots;\, n=0,1,2,\ldots\}$ one can
construct a standard Wiener process $\{W(t);\, t\geq 0\}$ with local time
process $\{\eta(x,t);\, x\in\mathbb R; \, t\geq 0\}$ such that, as
$n\to\infty$, we have for any $\varepsilon>0$
$$|S(n)-W(n)|=O(n^{1/4+\varepsilon})\quad {a.s.}
$$
and
$$
\sup_{x\in\mathbb Z}|\xi(x,n)-\eta(x,n)|=O(n^{1/4+\varepsilon})
\quad {a.s.},
$$
simultaneously.}

The following result about the 
uniformity of the local time is in Heyde \cite{H}, see also in Cs\'aki and R\'ev\'esz \cite{CR83}.

\medskip
\noindent
 {\bf Lemma D} (\cite{CR83}, \cite{H}) {\it
For the simple symmetric walk for any $\varepsilon>0$ we have}

$$\lim_{n\to \infty}
\frac{\sup_{x}|\xi(x+1,n)-\xi(x,n)|}{n^{1/4+\varepsilon}}=0\quad
a.s.
$$

\noindent
{\bf Remark 2.1} In fact \cite{CR83}  deals with more general random walks, but 
we only need it for a simple symmetric random walk.

The following result is the so called exponential Kolmogorov inequality. It is  a direct consequence of Doob's maximal inequality. Its proof can be found e.g. on page 139 of Williams \cite{W}.

\medskip\noindent
{\bf Lemma E}  Let {\it $X_j,\, j\geq 1$ be i.i.d. random variables with }
$\displaystyle{E\left(exp(\theta |X_ j|)\right)< \infty}$ for some $\theta>0$ and $E(X_j)=0.$ Then for any $\lambda>0$

$${\bp} \left(\max_{1\leq j\leq n}
\left |\sum_{i=1}^j X_i \right|>\lambda\right)
\leq
\exp \left (-\lambda\theta\right)
\left(
 E\left(exp(\theta |X_ j|)^n\right) 
+E\left(exp(-\theta |X_ j|)^n\right)
 \right).
$$

\medskip
The following result is a generalization of an inequality of T\'oth \cite{TO}.

\medskip\noindent
 \begin{lemma} 
 {\it Let} $G_i$, $i=1,2,\ldots$ {\it be 
i.i.d. random variables with the common geometric distribution }
${\bf P}(G_i=k)=\alpha (1-\alpha)^k,\quad k=0,1,2...$ for some $0<\alpha<1.$ {\it Then for $n$ big enough}

$$ {\bp} \left(\max_{1\leq j\leq n}
\left |\sum_{i=1}^j \left(G_i-\frac{1-\alpha}{\alpha}\right)\right|>\lambda\right)\leq
2\exp\left (-\frac{\lambda^2\alpha^2}{{4(1-\alpha)n}}\right)
$$
for $0<\lambda<n a$ with some $a>0.$
\end{lemma}
{\bf Proof.}
 The common moment generating function  of $G_i$-s for $ i=1,2,\dots$ is 
 
$$ {\bf E} \left(e^{\theta G_i} \right)=\frac{\alpha}{1-e^{\theta}(1-\alpha )}$$
 \noindent
and \,${\be}(G_i)=\frac{1-\alpha}{\alpha}.$
Then 
$${\bf E} \left(e^{\theta( G_i-\frac{1-\alpha}{\alpha} )} \right)
=\frac{\alpha}
{(1 -e^{\theta}(1-\alpha ) )e^{\theta(\frac {1-\alpha}{\alpha} )}}.$$
There exists a constant $ \theta_1>0,$ such that  for $0<\theta<\theta_1$ we have that
$$1 -e^{\theta}(1-\alpha )=1-(1-\alpha)(1+\theta +\frac{\theta^2}{2} +O(\theta^3)),$$
and
$$ e^{\theta(\frac {1-\alpha}{\alpha} )}=1+\frac{1-\alpha}{\alpha}\theta+\left(\frac{1-\alpha}{\alpha}\right)^2\frac{\theta^2}{2}+O(\theta^3), $$
 \noindent
where $\theta_1$ above might depend on $\alpha.$
Now with elementary calculation we get, that

$$(1 -e^{\theta}(1-\alpha ) ) e^{\theta\left(\frac {1-\alpha}{\alpha} \right)}=\alpha - \frac{1-\alpha}{2 \alpha}\theta^2+O(\theta^3)=\alpha\left(1-\frac{1-\alpha}{2 \alpha^2}\theta^2\right)+O(\theta^3).$$
and
$${\bf E} \left(e^{\theta( G_i-\frac{1-\alpha}{\alpha} )} \right)=1+\theta^2 \frac{1-\alpha}{2 \alpha^2}+O(\theta^3)=
e^{ \frac{(1-\alpha)\theta^2 }{2 \alpha^2} }+O(\theta^3)\leq e^{ \frac{(1-\alpha)\theta^2 }{ \alpha^2}}, $$
again with $0<\theta\leq\theta_2\leq\theta_1.$

Considering now the moment generating function  $\displaystyle{{\bf E} \left(e^{-\theta( G_i-\frac{1-\alpha}{\alpha})} \right)} $ and almost identical calculation results that
$${\bf E} \left(e^{-\theta( G_i-\frac{1-\alpha}{\alpha})} \right)\leq e^{ \frac{(1-\alpha)\theta^2 }{ \alpha^2}}, $$
as well, again with $0<\theta\leq\theta_3\leq\theta_2.$

Applying now the exponential Kolmogorov Inequality (Lemma E) we get that

\begin{eqnarray} && {\bf P} \left(\max_{1\leq j\leq n} \left | \sum_{i=1}^j \left( G_i- \frac{1-\alpha}{\alpha}\right)  \right |>\lambda \right)  \nonumber \\
 < && e^{-\lambda \theta} \left( {\bf E} \left(e^{\theta( G_i-\frac{1-\alpha}{\alpha})} \right)^n+{\bf E} \left(e^{-\theta( G_i-\frac{1-\alpha}{\alpha})} \right)^n \right) \nonumber \\
&&\leq2e^{-\lambda \theta} e^{n \frac{(1-\alpha)\theta^2 }{ \alpha^2}}\leq  2 e^{-\frac{\lambda^2 \alpha^2}{4n(1-\alpha) }}
\end{eqnarray} 
where we selected $\displaystyle{\theta=\frac{\lambda \alpha^2}{2n(1-\alpha)} }$, which obviously can be done for 
$n$ big enough. $\Box $

\smallskip

\section{Proof of the  Theorem 1.1}
\renewcommand{\thesection}{\arabic{section}} \setcounter{equation}{0}
\setcounter{thm}{0} \setcounter{lemma}{0}

First we are to redefine our random walk $\{{\mathbf C}(N);\,
N=0,1,2,\ldots\}$. It will be seen that the process described right below
is equivalent to that given in the Introduction.
For notational convenience we will say that the set $B$ contains the levels $\{m_j,\quad j=1,2, \ldots K \}.$
To begin with, on a suitable probability space consider two independent
simple symmetric (one-dimensional) random walks $S_1(\cdot)$, and
$S_2(\cdot)$. We may assume that on the same probability space we have 
$K$ sequences of i.i.d. geometric random variables 
$\{G_i(m_j),\,i\geq 1,\quad m=1,2, \ldots K \}$ which are independent from each other and $S_1(\cdot)$, and
$S_2(\cdot),$ with
$$
\mathbf{P}(G_i(m_j)=k)=\alpha_j (1-\alpha_j)^{k},\,\, k=0,1,2,\ldots \quad{\rm with} \quad 0<\alpha_j<1, \quad
   j=1,2, \ldots K\quad \label{geoseq}
$$
For simplicity we will say that the sequence of geometric random variables $\{G(m_j)_i,\,i\geq 1\}$ belongs to the horizontal level $m_j.$
We now construct our walk $\mathbf{C}(N)$ as follows. We will take all
the horizontal steps consecutively from $S_1(\cdot)$ and all the vertical
steps consecutively from $S_2(\cdot).$ 
Consider our walk starting from the origin. First the walk moves  vertically
until it arrives at a level which belongs to $B.$ (It is possible that no vertical step is needed, as the $x$-axis might 
belong to B).  If this level is  level $m_j\in B,$ then  it takes $G_1(m_j) $ horizontal steps
 from $S_1(\cdot).$ (Note that $G_1(m_j)=0$ is possible with probability $ \alpha_j$). Then we again take vertical  steps   from $S_2(\cdot),$ as needed to get to a level belonging to 
$B,$ then again some horizontal steps from $S_1(\cdot)$ as follows. If this level in $B$ is the  same as previously
then we take $G_2(m_j) $ horizontal steps. However if it is another level, lets say $m_{\ell}\in B$ then it takes  $G_1(m_{\ell}) $
horizontal steps from  $S_1(\cdot),$ then  the walk moves vertically taking steps from $S_2(\cdot)$ until it hits again a level in  $B$, when again it moves horizontally taking steps from $S_1(\cdot).$
In general, whenever the walk arrives at the level in $ B$  
then it takes some horizontal steps, the number of which is given by the next in 
line (first unused) geometric random variables  belonging to that level.

Let now $H_N,\, V_N$ be the number of horizontal and vertical steps, 
respectively, from the first $N$ steps of the just described process. 
Consequently, $H_N+V_N=N$, and
$$
\left\{{\bf C}(N);\, N=0,1,2,\ldots\right\}=
\left\{(C_1(N),C_2(N));\, N=0,1,2,\ldots\right\}
$$
\begg
\stackrel{{d}}{=}\left\{(S_1(H_N),S_2(V_N));\, N=0,1,2,\ldots\right\},
\label{equ}
\endd
where $\stackrel{{d}}{=}$ stands for equality in distribution.

Now we introduce a few more notations. Let $\xi_2(\cdot,\cdot)$ denote the 
local time of $S_2(\cdot)$.
\begin{eqnarray}
H_N&=&\#\{ k: \quad 1 \leq k\leq N, \,\,\, C_1(k)\neq C_1(k-1)\}\\
V_N&=& \#\{ k: \quad 1 \leq  k\leq N, \,\,\, C_2(k)\neq C_2(k-1)\}\\
D_2(V_N)&=&\sum_{j=1}^K \xi_2(m_j,V_N) \frac{1-\alpha_j}{\alpha_j} . \label{dset}
\end{eqnarray}

\noindent  Clearly  $H_N$
and $V_N$ are the number of horizontal and vertical steps, respectively, in 
the first $N$ steps of ${\bf C}(\cdot).$  $D_2(V_N)$ is the expected occupation time of
the levels belonging to  $B$ by ${\bf C}(\cdot)$ in the first $N$ steps.

\begin{lemma} For any $\ep>0$, as $N\to\infty$,
$$\max_{1\leq i\leq N} |H_i-D_2(V_i) |=O (N^{1/4+\ep})\quad a.s.$$
\end{lemma}

\noindent

\noindent
{\bf Proof.}  Recall that $H_N$ is the number of horizontal steps in our 
construction, and horizontal steps only occur on levels belonging to $B.$ When 
the vertical walk arrives to such a level, $m_j$, it takes some horizontal steps, 
the number of which follows geometric distribution with expected value  $\frac{1-\alpha_j}{\alpha_j}.$ 
In $V_N$ steps  the vertical walk $S_2(\cdot)$ spends $\xi_2(m_j,V_N)$ steps on the level $m_j$, thus the number of horizontal steps on this level  is the sum of $\xi_2(m_j,V_N)$ geometric random variables with  common expected value 
$\frac{1-\alpha_j}{\alpha_j}.$ 
The total number of horizontal steps is 
$\sum_{j=1}^K\sum_{i=1}^{\xi_2(m_j,V_{\ell})}G_i(m_j).$
However this statement is 
slightly incorrect, as if the $N$-th step is a horizontal one, the 
corresponding last geometric random variable might remain truncated. Denote by 
$H_N^+$ the number of horizontal steps which includes all the steps of this 
last geometric random variable. Then 
  $$H_N^+ -D_2(V_N)=\sum_{j=1}^K\sum_{i=1}^{\xi_2(m_j,V_N)}\left(G_i(m_j)-\frac{1-\alpha_j}{\alpha_j}\right), $$
where $G_i(m_j), \, i=1,2...$ are the i.i.d. geometric random variables, belonging to level $m_j,$ as
in Lemma 2.1. According to this lemma and Lemma A we have
\begin{eqnarray}&&
{\bf P}\left(\max_{1\leq \ell\leq N} |H_{\ell}^+ -D_2(V_{\ell})|>\lambda \right)=
{\bf P}\left(\max_{1\leq \ell\leq N} \left|\sum_{j=1}^K\sum_{i=1}^{\xi_2(m_j,V_{\ell})}\left(G_i(m_j)-\frac{1-\alpha_j}{\alpha_j}\right)\right|> \lambda \right)  \nonumber \\
&&< \sum_{j=1}^K {\bf P}\ \left(\max_{1\leq \ell\leq N}\left|\sum_{i=1}^{\xi_2(m_j,V_{\ell})}\left(G_i(m_j)-\frac{1-\alpha_j}{\alpha_j}\right)\right|>\frac{\lambda}{K}\right)  \nonumber \\
&& \leq\sum_{j=1}^K 2 \exp\left (-\frac{\lambda^2\alpha_j^2}{{4K^2(1-\alpha_j)N^{1/2+\epsilon}}}\right)
 \leq 2K \exp \left(-\frac{\lambda^2(\alpha^*)^2}{{4K^2(1-\alpha^*)N^{1/2+\epsilon}}}\right),\
 \end{eqnarray}
for any $\epsilon>0,$ with $\alpha^*=\min_{1 \leq j\leq K} \alpha_j$  as  the function  $\frac{\alpha^2}{1-\alpha}$ is increasing for $0<\alpha<1$, where we used the fact that 
$$\max_{1\leq \ell \leq N} \xi_2(m_j,V_{\ell}) \leq N^{1/2+\ep} \qquad a.s.$$ if $N$ is big enough.

Selecting $\lambda=N^{1/4+\ep}$ we get by the Borel-Cantelli lemma that for $N $ large enough
$$\max_{1\leq \ell\leq N} |H_{\ell}^+ -D_2(V_{\ell})|=O( N^{1/4+\ep}) \quad a.s.$$

Now to estimate the difference of $H_N$ and $H_N^+$, we have to observe, that their difference is not more than one single geometric random variable, which happens to be the last one used up to $N$.  Thus with $\alpha^*$ defined above we have
$${\bf P} (|H_N^+ -H_N|\geq N^{\ep})\leq {\bf P} ( \max_{1\leq j\leq K} \max _{1\leq i \leq N }G_i(m_j)\geq N^\ep)\leq
N \max_{1\leq j\leq K} {\bf P} (G_1(m_j)\geq N^\ep)\leq N  (1-\alpha^*)^{N^\ep}
$$
and hence by the Borel-Cantelli lemma
$$
H_N^+-H_N\leq N^{\varepsilon}\quad a.s.
$$
for all large $N,$ proving our lemma.
 $\Box$

Now observe that based on Lemma D
$$D_2(V_N)=\sum_{j=1}^K \xi_2(m_j,V_N) \frac{1-\alpha_j}{\alpha_j}=\xi_2(0,V_N) \sum_{j=1}^ K \frac{1-\alpha_j}{\alpha_j}  +O(N^{1/4+\ep})\quad a.s. $$

 Recall  the notation  $\displaystyle{A_K= \sum_{j=1}^ K \frac{1-\alpha_j}{\alpha_j}}$ given  in (\ref{alfa}).
Using this,  we have

$$D_2(V_N)=A_K \xi_2(0,V_N)+O(N^{1/4+\ep})\quad a.s.$$ implying by Lemma 3.1 that 
\begg H_N=A_K \xi_2(0,V_N) +O(N^{1/4+\ep})\quad a.s.\label{haen}   \endd 
as well. Moreover by (\ref{haen}) and Lemma A  we have that
 $H_N=O(N^{1/2+\ep})$ \, a.s., thus
$$|V_N-N|=O(N^{1/2+\ep}) \quad a.s., $$
\noindent
implying by the last statement of Lemma B that 
$$H_N=A_K\xi_2(0,N)+O(N^{1/4+\ep})\quad a.s.$$
Then we have from (\ref{equ}), and Lemmas B  and C that
\begin{eqnarray}
C_1(N)&=&S_1(H_N)=S_1(A_K \xi_2(0,N)+O(N^{1/4+\ep})) \nonumber \\
&=&W_1(A_K \xi_2(0,N)+O(N^{1/4+\ep}) )+O(N^{1/8+\ep}) \nonumber \\
&=&W_1(A_K \eta_2(0,N)+O(N^{1/4+\ep}))+O(N^{1/8+\ep})  \nonumber \\
&=&W_1(A_K \eta_2(0,N))+O(N^{1/8+\ep}) \quad a.s.
\end{eqnarray}
and 
$$C_2(N)=S_2(V_N)=S_2(N+O(N^{1/2+\ep}))= W_2 (N)+O(N^{1/4+\ep}) \quad a.s.$$
proving  our theorem. $\Box$

\section{Consequences}
\renewcommand{\thesection}{\arabic{section}} \setcounter{equation}{0}
\setcounter{thm}{0} \setcounter{lemma}{0}

\noindent
Define the continuous version  of our random walk process  on ${\mathbb C^2_K}$  by linear interpolation, 
as follows:
$$\{{\bf C}(xN)=(C_1(xN),C_2(xN)): 0\leq x\leq1\}.$$  
We have almost surely, as $N\to\infty$,
$$ \sup_{0\leq x\leq 1}\left\Vert\left(\frac{C_1(xN)-W_1(A_K\eta_2(0,xN))}
{N^{1/4}(\log\log N)^{3/4}},
\frac{C_2(xN)-W_2(xN)}{(N\log \log N)^{1/2}}\right)\right\Vert\to 0.$$

We have the following laws of the iterated logarithm (for the first statement see 
Theorem 2.2 in Cs\'aki {\it et al.}  \cite{CCFR95}).

$$\limsup_{n\to \infty} \frac{C_1(N)}{\sqrt{A_K}N^{1/4}(\log\log N)^{3/4}}
=\frac{2^{5/4}}{3^{3/4}}  \quad a.s. \quad {\rm and} \quad
\limsup_{N\to \infty} \frac{C_2(N)}{(2N\log\log N)^{1/2}}=1 \quad a.s.$$

As to the liminf behavior of the max functionals of
the two components,  we have the same results as for the two dimensional comb 
lattice \cite{CCFR08}. These results are based on the corresponding ones  
for Wiener process  and  the iterated process $W_1(\eta_2(0,t))$ and  
the work of Chung \cite {CH}, Hirsch \cite{HI}, Bertoin \cite{BER}, and 
Nane \cite{NA} . 

Based on \cite{NA}, we get the following: Let $\rho(n),\, n=1,2,\ldots$, be 
a non-increasing sequence of positive numbers such that $n^{1/4}\rho(n)$ is
non-decreasing. Then we have almost surely that
$$
\liminf_{N\to\infty}\frac{\max_{0\leq k\leq
N}C_1(k)}{N^{1/4}\rho(N)}=0\quad or\quad \infty
$$
and
$$
\liminf_{n\to\infty}\frac{\max_{0\leq k\leq
N}C_2(k)}{N^{1/2}\rho(N)}=0\quad or\quad \infty,
$$
according as to whether the series $\sum_1^\infty \rho(n)/n$ diverges or
converges.

\begin{equation}
\liminf_{N\to\infty}\left(\frac{8\log\log N}{\pi^2 N}\right)^{1/2}
\max_{0\leq k\leq N}|C_2(k)|=1 \qquad a.s.
\end{equation}

On the other hand, for the max functional of $|C_1(\cdot)|$  we obtain 
from \cite{CCFR08} the following result.

Let $\rho(n),\, n=1,2,\ldots$, be a non-increasing
sequence of positive numbers such that $n^{1/4}\rho(n)$ is
non-decreasing. Then we have almost surely that
$$
\liminf_{N\to\infty}\frac{\max_{0\leq k\leq
N}|C_1(k)|}{N^{1/4}\rho(N)}=0\quad or\quad \infty,
$$
as to whether the series $\sum_{n=1}^\infty \rho^2(n)/n$ diverges 
or converges.

\end{document}